\documentclass{article}
\usepackage[utf8]{inputenc}
\usepackage[T1]{fontenc}
\usepackage{mathtools, amsfonts, amssymb, amsthm}
\usepackage[maxbibnames=99]{biblatex}
\newtheorem{theorem}{Theorem}
\newtheorem{definition}{Definition}
\newtheorem{proposition}{Proposition}
\newtheorem{lemma}{Lemma}
\newtheorem{corollary}{Corollary}
\newtheorem{example}{Example}
\newtheorem{remark}{Remark}
\bibliography{biblio.bib}

\title{Laplacian operator on statistical manifold}
\author{Ruichao Jiang, Javad Tavakoli, Yiqiang Zhao}
\date{Feb 14, 2022}
\begin{document}

\maketitle
\begin{abstract}
    In this paper, we define a Laplacian operator on a statistical manifold, called the vector Laplacian. This vector Laplacian incorporates information from the Amari-Chentsov tensor. We derive a formula for the vector Laplacian. We also give two applications using the heat kernel associated with the vector Laplacian.
\end{abstract}
\section{Introduction}
    Information geometry studies statistical manifolds. In contrast to Riemannian geometry, a pair of dual affine connections $(\nabla, \Tilde{\nabla})$ is defined on a statistical manifold $M$. This dualistic structure can be described equivalently by a third order tensor $C=\nabla g$, where $\nabla$ is the covariant derivative operator. The $C$ is called the Amari-Chentsov tensor in information geometry and the non-metricity tensor in physics. The geometric meaning of $C$ is related to the change of length scale during an infinitesimal parallel transport.
    
    In this paper, we define a vector Laplacian operator on a statistical manifold with volume density. We explicitly consider the volume density for three reasons:
    \begin{enumerate}
        \item To apply information geometry on Bayesian statistics \cite{jiang}, one needs to specify a prior distribution. The choice of the prior distribution reflects the prior knowledge before inference. When one has completely no knowledge about the parameters, the Riemannian volume density (Jeffreys prior) is often used.
        \item Laplacian operators are usually $L^2$ symmetric operators:
        \begin{equation*}
            \langle f,\Delta g\rangle=\langle\Delta f,g\rangle,
        \end{equation*}
        where $f,g$ are functions and $\langle\cdot,\cdot\rangle$ is an inner product. This property is crucial for the semigroup theory and the existence of the heat kernel associated with $\Delta$. A volume density or a measure can be used to define the inner product above.
        \item In Riemannian geometry, the Riemannian volume form can be characterized as the unique volume form parallel with respect to the Levi-Civita connection. However, a non-Levi-Civita connection $\nabla$ does not necessarily have a parallel volume form. In fact, $\nabla$ has a parallel volume form if and only if its Ricci curvature is symmetric \cite{nomizu_book}.
    \end{enumerate}
    When defining the vector Laplacian $\Delta$, we follow two goals. First, $\Delta:TM\to TM$ should be an $L^2$ symmetric operator on the tangent bundle $TM$. Second, we want to incorporate the dualistic structure on $M$ into $\Delta$. As a result, we define $\Delta$ as the Bochner Laplacian on $M$. We then derive a formula for $\Delta$. In fact, our approach works for any vector bundle $E\to M$, not necessarily the tangent bundle. 
    
    In recent years, the heat kernel associated with some Laplacian operator plays an important role in applied areas such as example computer graphics \cite{sharp}, computer vision \cite{batard, singer}, dimensionality reduction \cite{coifman,}, and natural language processing \cite{lafferty1, lafferty2}. In those applications, the heat kernel $p_t(x,y)$ associated with the Laplacian operator can be used to define a distance function, different from the geodesic distance function. Assuming the compactness of $M$, the heat kernel associated with our vector Laplacian exists. Hence we can use $p_t(x,y)$ to define a distance function which incorporates the dualistic structure on $M$. Note that it is difficult to incorporate the dualistic structure in the standard approach to define the geodesic distance because $\nabla$ is not metric-compatible. As a result, the length of a $\nabla$-geodesic is not well-defined.
    
    The organization of the paper is as follows: Section~2 is an introduction to information geometry; in Section~3 we define the Laplacian operator $\Delta$ in information geometry and derive a formula for it; in Section~4 we discuss application; and we discuss our compactness assumption in Section~5.
\section{Information geometry and Weyl geometry}
In this section, we review some facts of information geometry \cite{amari, ay, nielson}.  In fact, a statistical manifold can be defined completely in an abstract way. However, to show its relevance to statistics, we adopt here a more concrete way: we start with a statistical model and work in local coordinate to define geometric objects in information geometry.

Consider a statistical model $\mathcal{P}$, which is a set of parametric densities $\mathcal{P}=\{p(\cdot|\theta)\}.$ We denote by $\Theta$ the parameter space, which is the set of all possible values that $\theta$ can take. $\Theta$ can be made into a Riemannian manifold by the following metric.
\begin{definition}
    In local coordinates, the Fisher information metric tensor $g_{ij}$ is defined by
    \begin{equation}
        g_{ij}=\text{E}_{\theta}\left[\partial_{i}\ell\partial_{j}\ell\right]
    \end{equation}
    where $E_{\theta}$ is the the transition kernel $X\times\Theta\to{[0,\infty)}$, $\ell$ is the log likelihood function, and $\partial_{i}$ is the partial derivative with respect to coordinate $i$.
\end{definition}
In addition to the Riemannian metric, the following tensor is defined in information geometry.
\begin{definition}
    In local coordinates, the Amari-Chentsov tensor $C_{ijk}$ is defined by
    \begin{equation}
    C_{ijk}=\text{E}_{\theta}\left[\partial_{i}l\partial_{j}l\partial_{k}l\right].
    \end{equation}
\end{definition}
From the Amari-Chentsov tensor $C$ and the Levi-Civita connection $\nabla^g$, a dual pair of affine connections on $\Theta$ can be defined.
\begin{definition}[Dual connections]
    Let $\nabla$ be an arbitrary torsion-free affine connection on a Riemannian manifold $(M,g)$. The connection $\Tilde{\nabla}$ dual to $\nabla$ is defined by the following equation:
    \begin{equation}
        \label{dual_metric_compatibility}
        Xg(Y,Z)=g(\nabla_XY,Z)+g(Y,\Tilde{\nabla}_{X}Z),
    \end{equation}
    where $p\in{M}$ and $X,Y,Z\in{T_{p}M.}$
\end{definition}
\begin{remark}
    The Levi-Civita connection $\nabla^g=\frac{\nabla+\Tilde{\nabla}}{2}$ is self-dual.
\end{remark}
The equivalence of the Amari-Chentsov tensor and the duality structure follows from the following proposition.
\begin{proposition}
    Let $(M,g,\nabla)$ be a statistical manifold and let $\Tilde{\nabla}$ be the dual connection of $\nabla$ with respect to $g$. Then, $C=\nabla g$. In local coordinates, let the Christoffel symbols of $\nabla$ and $\Tilde{\nabla}$ be $\Gamma_{ijk}$ and $\Tilde{\Gamma}_{ijk}$ respectively. Then
    \begin{equation*}
        \Tilde{\Gamma}_{ijk}-\Gamma_{ijk}=C_{ijk}.
    \end{equation*}
\end{proposition}
Following \cite{zhang}, we use the difference tensor $K$ to represent $C$ as an operator.
\begin{definition}
    \label{difference}
    Let $(M,g,\nabla)$ be a statistical manifold with dual connection $\Tilde{\nabla}$. The difference tensor $K$ is defined as follows:
    \begin{equation*}
        K_XY=\Tilde{\nabla}_XY-\nabla_XY.
    \end{equation*}
\end{definition}
The relation between $K$ and $C$ is given by the following equation:
\begin{equation}
    C(X,Y,Z)=g(K_XY,Z).
\end{equation}
\section{Vector Laplacian in information geometry}
    In this section, we define a Laplacian operator $\Delta$ acting on vector fields on the statistical manifold $(M,g,\nabla)$ with density $\rho=\text{e}^{-f}\sqrt{\det{g}}$. We impose the following condition: $\Delta$ should be an $L^2$ symmetric operator with respect to the volume form $\rho$, i.e.
    \begin{equation}
        \label{symmetric}
        \int_Mg(\Delta X,Y)\rho=\int_Mg(X,\Delta Y)\rho
    \end{equation}
    This requirement is crucial because being symmetric is necessary to apply semigroup theory to show the existence of the heat kernel.
    
    We recall the definition of the divergence operator $\text{div}_f$ on a weighted manifold.
    \begin{definition}
        \label{div}
        Let $(M,g,\rho)$ be a weighted manifold with density $\rho=\text{e}^{-f}\sqrt{\det{g}}$. The divergence operator $\text{div}_f$ is defined as the (negative) formal adjoint operator of the exterior derivative operator d with respect to $\rho$.
        \begin{equation*}
            \int_M\left(X,\text{d}h\right)\rho=-\int_M\text{div}_f(X)h\rho,
        \end{equation*}
        where $X$ is an arbitrary vector field, $h$ is a scalar function, and $(\cdot,\cdot)$ is the canonical pairing of vector field and $1$-form. In local coordinates, we have
        \begin{equation*}
            \text{div}_f(X)=\frac{1}{\rho}\partial_i(\rho X^i)=\frac{1}{\text{e}^{-f}\sqrt{\det{g}}}\partial _i(\text{e}^{-f}\sqrt{\det{g}}X^i).
        \end{equation*}
    \end{definition}
    \begin{remark}
        For a Riemannian manifold $(M,g)$ with the Levi-Civita connection $\nabla^g$ and the Riemannian volume form $\sqrt{\det{g}}$, the divergence operator can also be characterized by the covariant derivative
        $\text{div}(X)=\frac{1}{\sqrt{\det{g}}}\partial _i(\sqrt{\det{g}}X^i)=\nabla^g_iX^i$. However, such a coincidence does not hold for a general affine connection with torsion or Amari-Chentsov tensor.
    \end{remark}
    The divergence operator $\text{div}_f$ on a weighted manifold is defined only by the differentiable sturcture and the volume form, not by any affine connection structure. Hence, we use Definition \ref{div} as the definition of the divergence operator on a statistical manifold.
    We collect some formulas for the divergence operator on a statistical manifold.
    \begin{lemma}
        \label{divergence}
        The divergence operator $\text{div}_f$ satisfies the following formulas:
        \begin{enumerate}
            \item $\text{div}_f(X)=\text{div}(X)-Xf$;
            \item $div_f(hX)=h\text{div}_f(X)+Xh$;
            \item $\int_M\text{div}_f(X)\rho=0$ for all vector fields $X$.
        \end{enumerate}
    \end{lemma}
    \begin{proof}
        (1) We compute,
        \begin{equation*}
            \begin{split}
                \text{div}_f(X)&=\frac{1}{\text{e}^{-f}\sqrt{\det{g}}}\partial _i(\text{e}^{-f}\sqrt{\det{g}}X^i)\\
                &=\frac{1}{\sqrt{\det{g}}}\partial _i(\sqrt{\det{g}}X^i)-\partial_ifX^i\\
                &=\text{div}(X)-Xf.
            \end{split}
        \end{equation*}
        (2) Using $(1)$,
        \begin{equation*}
            \begin{split}
                \text{div}_f(hX)&=\text{div}(hX)-hXf\\
                &=h\text{div}(X)+Xh-hXf\\
                &=h\text{div}_f(X)+Xh.
            \end{split}
        \end{equation*}
        (3) Let $1$ be a constant function on $M$. Then,
        \begin{equation*}
            \int_M\text{div}_f(X)\rho=\int_M\text{div}_f(X)1\rho=\int_M\left(X,\text{d}1\right)\rho=0.
        \end{equation*}
    \end{proof}
    To satisfy the symmetric condition in equation (\ref{symmetric}), we use Bochner Laplacian operator $\Delta\coloneqq\nabla^*\nabla$, where $\nabla^*$ is the formal adjoint of the covariant derivative $\nabla:\Gamma\left(TM\right)\to\Gamma\left(T^*M\otimes TM\right)$, considered as a map between sections of bundles. The following theorem gives the expression of the adjoint connection on a statistical manifold.
    \begin{theorem}
        \label{adjoint}
        The adjoint of the covariant derivative is given by the following formula:
        \begin{equation*}
            \nabla^*_XZ=-\Tilde{\nabla}_XZ-\text{div}_f(X)Z,
        \end{equation*} 
        where $\Tilde{\nabla}$ is the dual connection.
    \end{theorem}
    \begin{proof}
         From equation ($\ref{dual_metric_compatibility}$) and Lemma \ref{divergence}, we obtain,
        \begin{equation*}
            g(\nabla_XY,Z)=\text{div}_f(g(Y,Z)X)-g(Y,Z)\text{div}_fX-g(Y,\Tilde{\nabla}_XZ).
        \end{equation*}
        Integrate the above equation to obtain
        \begin{equation*}
            \int_Mg(\nabla_XY,Z)\rho=\int_M\text{div}_f(g(Y,Z)X)\rho-\int_Mg\left(Y,(\text{div}_fX+\Tilde{\nabla}_X)Z\right)\rho.
        \end{equation*}
        The first term on the RHS is a total divergence term and vanishes by Lemma \ref{divergence}. Therefore,
        \begin{equation*}
            \int_Mg(\nabla_XY,Z)\rho=\int_Mg\left(Y,\left(-\text{div}_fX-\Tilde{\nabla}_X\right)Z\right)\rho.
        \end{equation*}
        Hence $\nabla^*_X=-\Tilde{\nabla}X-\text{div}_fX$.
    \end{proof}
    \begin{remark}
        The adjoint connection $\nabla^*:\Gamma(T^*M\otimes TM)\to\Gamma(TM)$ should act on the sections of the tensor product of the cotangent bundle and the tangent bundle (or more generally, some vector bundle). Hence, we can also write $\nabla^*(X^\flat\otimes Z)=\nabla^*_XZ=-\Tilde{\nabla}_XZ-\text{div}_f(X)Z$.
    \end{remark}
    With the adjoint connection, we define the vector Laplacian on a statistical manifold with density as follows.
    \begin{definition}
        \label{vector_laplacian_definition}
        Let $(M,g,\nabla,\rho)$ be a statistical manifold with density and let $\nabla^*$ be the formal adjoint of $\nabla$. The Laplacian operator is
        \begin{equation*}
            \Delta X\coloneqq\nabla^*\nabla X
        \end{equation*}
    \end{definition}
    We derive a formula for the Laplacian operator on the statistical manifold with density.
    \begin{theorem}
    \label{laplacian_formula}
    The vector Laplacian on a statistical manifold $(M,g,\nabla)$ with density $\rho=\text{e}^{-f}\sqrt{\det{g}}$ is given by
        \begin{equation*}
            X=-\text{Tr}\left(\text{Hess}\left(X\right)\right)-\frac{1}{2}K^i\nabla_iX+g^{ij}\partial_jf\nabla_iX.
        \end{equation*}
    \end{theorem}
    \begin{proof}
        We have
        \begin{equation*}
            \begin{split}
                &\nabla X=dx^i\otimes\nabla_iX\\
                &\nabla^*(Y^\flat\otimes X)=-\Tilde{\nabla}_YX-\text{div}_f(Y)X,
            \end{split}
        \end{equation*}
        where $X,Y$ are vector fields, and $\flat$ is the music isomorphism: $\left(\partial_j\right)^\flat=g_{ij}dx^i$.
        Then,
        \begin{equation*}
            \begin{split}
                \Delta X&=\nabla^*\left(dx^i\otimes\nabla_iX\right)\\
                &=\nabla^*\left[g^{ij}\left(\partial_j\right)^\flat\otimes\nabla_iX\right]\\
                &=\nabla^*\left[\left(\partial_j\right)^\flat\otimes g^{ij}\nabla_iX\right]\\
                &=-\Tilde{\nabla}_j\left(g^{ij}\nabla_iX\right)-\text{div}_f\left(\partial_j\right)g^{ij}\nabla_iX\\
                &=-\left(\partial_jg^{ij}+g^{ij}\Tilde{\nabla}_j+g^{ij}\text{div}_f(\partial_j)\right)\nabla_iX.
            \end{split}
        \end{equation*}
        By the above equation and Lemma \ref{adjoint}, we obtain
        \begin{equation}
            \label{vector_bundle}
            \begin{split}
                \Delta X&=-g^{ij}\Tilde{\nabla}_j\nabla_iX-\left(\partial_jg^{ij}+g^{ij}\partial_j\log{\sqrt{\det{g}}}\right)\nabla_iX+g^{ij}\partial_jf\nabla_iX\\
                &=-g^{ij}\Tilde{\nabla}_j\nabla_iX-\frac{1}{\sqrt{\det{g}}}\partial_j\left(\sqrt{\det{g}}g^{ij}\right)\nabla_iX+g^{ij}\partial_jf\nabla_iX.\\
            \end{split}
        \end{equation}
        Using the formula $g^{ij}\left(^k_{ij}\right)=-\frac{1}{\sqrt{\det{g}}}\partial_j\left(\sqrt{\det{g}}g^{jk}\right)$ from Riemannian geometry, where $\left(^k_{ij}\right)$ are the Christoffel symbols of the Levi-Civita connection, we obtain
        \begin{equation}
            \label{tangent_bundle}
            \begin{split}
                \Delta X&=-g^{ij}\Tilde{\nabla}_j\nabla_iX+g^{ij}\left(^k_{ij}\right)\nabla_kX+g^{ij}\partial_jf\nabla_iX\\
                &=-g^{ij}\left(\nabla_j\nabla_iX-\Gamma^k_{ij}\nabla_kX\right)-\frac{1}{2}g^{ij}K^k_{ij}\nabla_kX+g^{ij}\partial_jf\nabla_iX\\
                &=-\text{Tr}\left(\text{Hess}\left(X\right)\right)-\frac{1}{2}K^i\nabla_iX+g^{ij}\partial_jf\nabla_iX
            \end{split}
        \end{equation}
    \end{proof}
    \begin{remark}
        In fact, all the derivations up to equation (\ref{vector_bundle}) hold true for a general vector bundle if we use a bundle metric in equation (\ref{symmetric}) and pay attention to the summing indices when the rank of the bundle is different from the dimension of the manifold.
    \end{remark}
    \begin{remark}
        Equation (\ref{tangent_bundle}) does not say that the Amari-Chentsov tensor only appears as an additional term $\frac{1}{2}K^i\nabla_iX$ because Hess$(X)$ also depends on it.
    \end{remark}
    We immediately have the following corollary.
    \begin{corollary}
        If $\nabla$ is compatible with the metric and $f=0$, equation (\ref{vector_bundle}) says that Bochner Laplacian coincides with the connection Laplacian in Riemannian geometry.
    \end{corollary}
    \iffalse
    Since equation (\ref{vector_bundle}) holds true for a general vector bundle and a scalar function can be identified with a section of a trivial line bundle, our vector Laplacian also applies to functions. We have the following corollary.
    \begin{corollary}
        On a trivial line bundle $\mathbb{R}\times M$ with a flat connection $\nabla=d$, the vector Laplacian is given by the following:
        \begin{equation}
            \label{line_bundle_laplacian}
            \Delta h=\Delta^g_fh=\Delta^gh-\partial^if\partial_ih,
        \end{equation}
        where $\Delta^g_f$ is the weighted Laplace operator and $\Delta^g$ is the Beltrami-Laplace operator.
    \end{corollary}
    \begin{proof}
        Let $h$ be a smooth function. We have $\Tilde{\nabla}_j\nabla_ih=\partial_i\partial_jh$. Then,
        \begin{equation*}
            \begin{split}
                \Delta h&=-g^{ij}\partial_i\partial_jh-\frac{1}{\sqrt{\det{g}}}\partial_j\left(\sqrt{\det{g}}g^{ij}\right)\partial_ih+g^{ij}\partial_jf\partial_ih\\
                &=\Delta^gh+\partial^if\partial_ih=\Delta^g_fh.
            \end{split}
        \end{equation*}
    \end{proof}
    \begin{remark}
       Even we tried to incorporate Amari-Chentsov tensor into $\Delta$, $\Delta$ does not depend on it in scalar case.
    \end{remark}
    \fi
    
\section{Heat kernel and its applications}
    In this section, we use the heat kernel of the vector Laplacian operator $\Delta$ to construct some useful quantities.
    
    Formally, from the vector Laplacian operator $\Delta:\Gamma(TM)\to\Gamma(TM)$, the heat semigroup $\text{e}^{t\Delta}:\Gamma(TM)\to\Gamma(TM)$ can be constructed. The heat kernel $p_t(x,y):T_yM\to T_xM$ is then the integral kernel of $\text{e}^{t\Delta}$:
    \begin{equation}
        \left(\text{e}^{t\Delta}X\right)(x)=\int_Mp_t(x,y)X(y)\rho(y).
    \end{equation}
    From equation (\ref{tangent_bundle}), $\Delta$ is an elliptic operator. This is easy to see from the symbol of $\Delta$. From our definition, $\Delta$ is symmetric (formally self-adjoint). To ensure that the heat kernel associated with $\Delta$ exists, $\Delta$ needs to be essentially self-adjoint, i.e. it possesses a unique self-adjoint extension. If the statistical manifold $M$ is compact, from functional analysis, we know that being elliptic and symmetric implies essential self-adjointness. Hence assuming compactness, we can have the following examples. 
    
    The first example generalizes the vector diffusion distance introduced in \cite{singer}.
    \begin{example}
        Following \cite{singer}, suppose that $-\Delta$ has discrete positive eigenvalues $\lambda_0<\lambda_1\leq\lambda_2\leq\cdots$ with eigenvectors $\Delta X_n=-\lambda_nX_n$, counted with multiplicity. This is permitted by the compactness assumption. Then we have the expansion $p_t(x,y)=\sum e^{-\lambda_nt}X_n(x)X_n(y)^\flat.$
        The vector diffusion distance on a statistical manifold $M$ is defined as follows:
        \begin{equation}
            \begin{split}
                d_t^2(x,y)&=\text{Tr}\left[p_t(x,x)\otimes p_t(x,x)+p_t(y,y)\otimes p_t(y,y)-2p_t(x,y)\otimes p_t(x,y)\right]\\
                &=\sum_{n=0}^\infty\sum_{m=0}^\infty e^{-\left(\lambda_n+\lambda_m\right)t}\left[g_x(X_n,X_m)-g_y(X_n,X_m)\right]^2.
            \end{split}
        \end{equation}
        That $d_t(x,y)$ is indeed a distance function follows from reproducing Hilbert space theory. $d_t(x,y)$ incorporates all information-geometric data, i.e. $g$, $\nabla$, and $\rho$ on $M$. Since $\nabla$ is not metric compatible, the length of a connection geodesic is not well-defined. Therefore, it is difficult to come up with a distance function that incorporates the affine connection data. The vector diffusion distance $d_t(x,y)$ associated with our vector Laplacian resolves this problem.
    \end{example}
    The next example is to deal with the situation where we have a parametric model $X\to\Theta$ from a sample space to a parameter space. This example generalizes the approach of \cite{lafferty1} to vector case.
    \begin{example}
    The parameter space $\Theta$ is equipped with a statistical manifold structure. A function on the product of the sample space $K:X\times X\to\mathbb{R}$ is called a kernel function if it is symmetric and positive-definite. The heat kernel on $\Theta$ can be used to define a kernel function on $X\times X$ as follows:
    \begin{equation}
        \label{kernel_function}
        \begin{split}
            {K}_t(x,x')&=\int_{\Theta\times\Theta}g\left(\nabla_\theta f(\theta|x),p_t(\theta,\theta')\nabla_{\theta'} f(\theta'|x')\right)\rho(\theta)\otimes\rho(\theta')\\
            &=\int_\Theta g\left( \nabla_\theta f(\theta|x),\text{e}^{t\Delta }\nabla_\theta f(\theta|x')\right)\rho(\theta),
        \end{split}
    \end{equation}
    where $f(\theta|x)=\frac{f(x|\theta)g(\theta)}{\int_\Theta f(x|\theta)\rho(\theta)}$ is the Bayesian posterior and by abuse of notation $\nabla_\theta f(\theta|x)$ denotes a vector field, not a $1$-form. It is easy to show that $K_t(x,x')$ is symmetric and positive-definite. Hence, $K_t(x,x')$ can be further used to define a distance function on $X$ as follows:
    \begin{equation}
        d_t^2(x,y)=K_t(x,x)+K_t(y,y)-2K_t(x,y).
    \end{equation}
    First of all, $K_t(x,x')$ includes all information geometric data in the parameter space, i.e. $g$, $\nabla$ and $\rho$. Moreover, the authors of \cite{lafferty1} actually used two prior distributions: they used the Laplace-Beltrami operator $\Delta^g$ and in particular they used the $L^2$ symmetric property of $\Delta^g$, which means that the Riemannian density $\sqrt{\det{g}}$ is used on $\Theta$ (otherwise their Laplacian would not be symmetric). Then, they used an arbitrary volume density $\rho$ to define the kernel function. Our $K_t(x,x')$ only depends on $\rho$.
    \end{example}
\section{Discussion}
One major assumption that we made is the compactness of the statistical manifold. However, not all statistical manifolds are compact. For example, the parameter space of Gaussian distributions is the Lobachevsky plane, which is non-compact. In Riemannian geometry, the next thing to consider after compactness is the geodesic completeness: the Laplacian-Beltrami operator is essentially self-adjoint on a geodesically complete Riemannian manifold \cite{yau}. This result also holds true for weighted manifolds \cite{grigor}. However, for a statistical manifold, there exist at least two notions of completeness: (1) Metric completeness, which means that the Levi-Civita connection $\nabla^g$ is complete. It is also equivalent to the metric space completeness with respect to the geodesic length function. (2) Affine connection completeness, which means that the solution of the geodesic equation for $\nabla$ can be extended to $(-\infty,+\infty)$ in time. Nomizu \cite{nomizu} found a manifold whose metric is incomplete but affine connection is complete. Furthermore, there exists a compact statistical manifold whose affine connection is not complete \cite{opozda}, which complicates the relevance between the completeness of $\nabla$ and the essential self-adjointness of $\Delta$.
\section{Conclusion}
In Riemannian geometry, various Laplacian operators encode geometric information in their heat kernels and spectra. In this paper we introduced the vector Laplacian operator $\Delta$ on a statistical manifold with density $(M,g,\nabla,\rho)$. To incorporate the contribution from the Amari-Chentsov tensor and ensure that $\Delta$ is symmetric, we defined $\Delta$ as the Bochner Laplacian on $M$. As stated in Section~1, there are three reasons for why we consider the density $\rho$ as an independent piece of data. We derived formulas for $\Delta$ and mentioned two examples which extend the existing applications of the heat kernel. We discussed the assumption that we made about the existence of the heat kernel, which shows the complexity of the problem at hand.
\printbibliography
\end{document}